\documentclass[11pt]{article}
\usepackage{amsmath}
\usepackage{amssymb}
\usepackage{amsfonts}

\makeatletter

\@addtoreset{equation}{section}
\makeatother
\def\<{\langle}
\def\>{\rangle}

\newtheorem{lem}{Lemma}[section]
\newtheorem{theo}{Theorem}[section]
\newtheorem{rem}{Remark}[section]

\makeatletter
   
   \@addtoreset{equation}{section}
\makeatother

\begin{document}
\title{\bf Asymptotic Profile of Solutions for Some Wave\\ Equations with Very Strong Structural Damping}
\author{Ryo IKEHATA\thanks{Corresponding author: ikehatar@hiroshima-u.ac.jp} and Shin IYOTA\\ {\small Department of Mathematics, Graduate School of Education, Hiroshima University} \\ {\small Higashi-Hiroshima 739-8524, Japan}}
\maketitle
\begin{abstract}
We consider the Cauchy problem in ${\bf R}^{n}$ for some types of damped wave equations. We derive asymptotic profiles of solutions with weighted $L^{1,1}({\bf R}^{n})$ initial data by employing a simple method introduced in \cite{Ik-3,Ik-4}. The obtained results will include regularity loss type estimates, which are essentially new in this kind of equation.
\end{abstract}

\section{Introduction}
\footnote[0]{Keywords and Phrases: Wave equation; Structural damping; Asymptotic profiles; Regularity loss; Fourier Analysis; Low frequency; High frequency; Weighted $L^{1}$-initial data.}
\footnote[0]{2010 Mathematics Subject Classification. Primary 35L15, 35L05; Secondary 35B40, 35B65.}

We are concerned with the Cauchy problem for wave equations in ${\bf R}^{n}$ ($n \geq 1$) with the structural damping term
\begin{equation}
u_{tt}(t,x) - \Delta u(t,x)  +(-\Delta)^{\theta}u_{t}(t,x) = 0,\ \ \ (t,x)\in (0,\infty)\times {\bf R}^{n} ,\label{eqn}
\end{equation}
\begin{equation}
u(0,x)= u_{0}(x), \quad u_{t}(0,x) = u_{1}(x), \quad x\in {\bf R}^{n},\label{initial}
\end{equation}
where $\theta > 1$. The initial data $u_{0}$ and $u_{1}$ are also chosen from the usual energy space (for simplicity)
\[[u_{0},u_{1}] \in H^{1}({\bf R}^{n}) \times L^{2}({\bf R}^{n}).\]
\noindent
This model equation (1.1) with $\theta \in [1,2]$ is recently introduced in the paper written by Ghisi-Gobbino-Haraux \cite{HGG} to study the (unique) existence of the global in time solutions and its smoothig effect for $t > 0$. In this sense, we call the model (1.1) as the GGH-model for short. Due to \cite[Theorem 2.1]{HGG}, it is known that the problem (1.1)-(1.2) admits a unique mild solution in the class\\
\[u \in C([0,+\infty);H^{1}({\bf R}^{n})) \cap C^{1}([0,+\infty);L^{2}({\bf R}^{n}))\]
satisfying the smoothing property:
\[\frac{\partial^{m}}{\partial t^{m}}u \in C((0,\infty);H^{1+2m(\theta-1)}({\bf R}^{n})) \quad (\forall m \in {\bf N}).\]
\noindent
As for the Cauchy problem of the structurally damped wave equation (1.1) with $\theta \in [0,1]$ many interesting results about the asymptotic behavior of solutions are well studied by many mathematicians.

(1)\,when $\theta = 0$, it is known that $u(t,x) \sim MG(t,x)$ or $u(t,x) \sim v(t,x)$ as $t \to \infty$ with some constant $M \ne 0$, where
\[G(t,x) = \displaystyle{\frac{1}{(\sqrt{4\pi t})^{n}}}e^{-\frac{\vert x\vert^{2}}{4t}},\]
and $v(t,x)$ is the (unique) solution to the heat equation
\[
v_{t}(t,x) - \Delta v(t,x) = 0,\ \ \ (t,x)\in (0,\infty)\times {\bf R}^{n}, 
\]
\[
v(0,x)= u_{0}(x) + u_{1}(x), \quad x\in {\bf R}^{n}.
\]
These results called as the diffusion phenomena are well-studied by Chill-Haraux \cite{CH}, Han-Milani \cite{HM}, Hayashi-Kaikina-Naumkin \cite{HKN}, Hosono \cite{H}, Hosono-Ogawa \cite{HO}, Ikehata-Nishihara \cite{IN}, Kawakami-Takeda \cite{KT}, Kawakami-Ueda \cite{KU}, Narazaki \cite{Na}, Nishiyama \cite{Nishiyama}, Radu-Todorova-Yordanov \cite{RTY}, Said-Houari \cite{said}, Takeda \cite{T}and Wakasugi \cite{W}, and the references therein. In particular, we should cite a deep result due to Nishihara \cite{N-2} such that
\[u(t,x) \sim v(t,x) + e^{-t/2}w(t,x)\quad (t \to \infty),\]
where $w(t,x)$ is the corresponding solution to the free wave equation  
\[
w_{tt}(t,x) - \Delta w(t,x) = 0,\ \ \ (t,x)\in (0,\infty)\times {\bf R}^{n}, 
\]
\[
w(0,x)= u_{0}(x), \quad w_{t}(0,x) = u_{1}(x), \quad x\in {\bf R}^{n}.
\]

(2)\, In the case of $\theta \in [0,1/2)$ Karch \cite{K} derived the asymptotic self-similar profile as $t \to \infty$ of the solutions. In fact, Karch also treated the nonlinear problems. Also, D'Abbicco-Reissig \cite{DR}, Lu-Reissig \cite{LR} and the references therein studied the decay estimates of various norms of solutions and the corresponding nonlinear problems in the case of general $\theta \in (0,1]$.

(3)\,In the case of $\theta \in [1/2,1]$, Ikehata-Natsume \cite{INatsume} and Char\~ao-da Luz-Ikehata\cite{CLI} derived the total energy and $L^{2}$ decay estimates of solutions to (1.1)-(1.2). Especially, for the Cauchy problem of the viscoelastic equation (1.1) with $\theta = 1$, recently Ikehata-Todorova-Yordanov \cite{ITY} in an abstract framework, and Ikehata \cite{Ik-4} in a concrete setting have derived its asymptotic profile such that
\[u(t,x) \sim  e^{-tA/2}(\cos(tA^{1/2})u_{0} + A^{-1/2}\sin(tA^{1/2})u_{1}), \quad (t \to +\infty),\]
where $A = -\Delta$ in $L^{2}({\bf R}^{n})$. This implies an oscillation property of the solution to (1.1) with $\theta = 1$. In this connection, Ponce \cite{p} and Shibata \cite{Shibata} derived the $L^{p}$-$L^{q}$ estimates of solutions to the equation (1.1) with $\theta = 1$ (viscoelastic equation case). 

From the observations (1)-(3) above, in the case of $\theta \in [0,1]$ one can say that more or less we have already known the asymptotic profile of the solution to problem (1.1)-(1.2).  But, these results are quite restricted to the lower power case such that $\theta \in [0,1]$. After GGH model with $\theta = 2$ is presented, it seems to be still open to discover the asymptotic profile of the solution to (1.1)-(1.2) with more general $\theta \in (1,2]$.\\    

Our main purpose is to find an asymptotic profile of the solution to problem (1.1)-(1.2) with $\theta = 2$, which is introduced by Ghisi-Gobbino-Haraux \cite{HGG}. Our results are read as follows.

\begin{theo}\,Let $n \geq 1$, $\theta > 1$ and $\ell \geq 0$. If $[u_{0},u_{1}] \in (H^{\ell+1}({\bf R}^{n}) \cap L^{1}({\bf R}^{n})) \times (H^{\ell}({\bf R}^{n}) \cap L^{1}({\bf R}^{n}))$, then it is true that
\[\Vert u_{t}(t,\cdot)\Vert^{2} + \Vert\nabla u(t,\cdot)\Vert^{2} \leq C(1+t)^{-\frac{n}{2\theta}}\Vert u_{1}\Vert_{1}^{2}+C(1+t)^{-\frac{n+2}{2\theta}}\Vert u_{0}\Vert_{1}^2\]
\[+C(1+t)^{-\frac{\ell}{\theta-1}}(\Vert D^{\ell}u_{1}\Vert^2 + \Vert D^{\ell+1} u_{0}\Vert^2),\]
where $C>0$ is a constant.
\end{theo}
\begin{theo}\,Let $n \geq 3$, $\theta > 1$ and $\ell \geq 1$. If $[u_{0},u_{1}] \in (H^{\ell}({\bf R}^{n}) \cap L^{1}({\bf R}^{n})) \times (H^{\ell-1}({\bf R}^{n}) \cap L^{1}({\bf R}^{n}))$, then it is true that
\[\Vert u(t,\cdot)\Vert^2 \leq C(1+t)^{-\frac{n-2}{2\theta}}\Vert u_{1}\Vert_{1}^{2}+C(1+t)^{-\frac{n}{2\theta}}\Vert u_{0}\Vert_{1}^2+C(1+t)^{-\frac{\ell}{\theta-1}}(\Vert D^{\ell-1}u_{1}\Vert^2 + \Vert D^{\ell} u_{0}\Vert^2),\]
where $C>0$ is a constant.
\end{theo}
Based on the observations above, we study the asymptotic profile of the solutions to problem (1.1)-(1.2) with $\theta = 2$ as an example and for simplicity. A generalization to $\theta > 1$ will be an easy exercise.   
\begin{theo}
Let $n\geq 1$, $\theta = 2$ and $\ell > \displaystyle{\frac{n}{4}-\frac{1}{2}}$ {\rm (}$n \geq 6${\rm )}, and $\ell \geq 1$ {\rm (}$1 \leq n \leq 5${\rm )}. If $[u_0,u_1]\in(H^{\ell}({\bf R}^n) \cap L^{1,1}({\bf R}^n))\times(H^{\ell-1}({\bf R}^n)\cap L^{1,1}({\bf R}^n))$, then the solution $u(t,x)$ to problem {\rm (1.1)-(1.2)} satisfies
\[\int_{{\bf R}^n}\vert\mathcal{F}(u(t,\cdot))(\xi)-\{P_1e^{-t\vert\xi\vert^4/2}\frac{\sin{(t\vert\xi\vert)}}{\vert\xi\vert}+P_0e^{-t\vert\xi\vert^4/2}\cos{(t\vert\xi\vert)}\}\vert^2d\xi\]
\[\leq C\Vert u_1\Vert_{1,1}^2 t^{-\frac{n}{4}} + C\Vert u_0\Vert_{1,1}^2t^{-\frac{n+2}{4}} + Ce^{-\alpha t}(\Vert u_1\Vert_1^2 + \Vert u_0\Vert_1^2+\Vert u_1\Vert^2+\Vert u_0\Vert^2)\]
\[+ Ct^{-\ell}(\Vert D^{\ell-1}u_{1}\Vert^2 + \Vert D^{\ell} u_{0}\Vert^2)\quad(t\gg1).\]
\end{theo}
\begin{rem}{\rm Totally speaking, the discovery of the regularity loss type structure in the high frequency region is completely new for this type of equations. The condition $\ell$ on the regularity of the initial data is based on the consideration from Lemmas 4.1 and 4.2 below, which imply that the right hand side of the inequality stated in Theorem 1.3 is the remainder term, and in this case the leading term of the solution $u(t,x)$ is the so called diffusion wave, i.e.,
\[{\cal F}^{-1}\left(P_1e^{-t\vert\xi\vert^4/2}\frac{\sin{(t\vert\xi\vert)}}{\vert\xi\vert}+P_0e^{-t\vert\xi\vert^4/2}\cos{(t\vert\xi\vert)}\right)(x).\]
$\ell \geq 1$ is necessary in order to guarantee the unique existence of mild solutions in the framework of the $H^{1}\times L^{2}$ initial data. We have another option how to choose a class of initial data  (see \cite{HGG}).
}
\end{rem}
\begin{rem}{\rm 
As a result, even in the case of $\theta = 2$, the asymptotic profile of the solution to problem (1.1)-(1.2) as $t \to \infty$ is almost the same as that of $\theta = 1$ derived in \cite{Ik-4}. However, we have just encountered the completely different aspects in the decay estimates shown in Theorems 1.1-1.3, compared with the case of $\theta \leq 1$, such that the regularity loss type estimates are appeared in those of the high frequency part of the remainder term. This is a big difference between $\theta \in [0,1]$ and $\theta \in (1,2]$. $\theta = 1$ is, in this sense, critical.}
\end{rem}
\begin{rem}{\rm 
From our results, we can present an open problem: can one find an asymptotic profile of the solution $u(t,x)$ to problem (1.1)-(1.2) with $\theta > 1$ in the case when the initial data have a low regularity?}
\end{rem}

Our plan in this paper is as follows. In section 2, we shall prove Theorems 1.1 and 1.2 by the energy method in the Fourier space due to \cite{UKS}, and in section 3 we prove Theorem 1.3 by the use of the method introduced by \cite{Ik-4}. As an application, we will discuss the optimality concerning the decay rate of the $L^{2}$-norm of solutions in Section 4.\\

{\bf Notation.}\,{\small Throughout this paper, $\| \cdot\|_q$ stands for the usual $L^q({\bf R}^{n})$-norm. For simplicity of notations, in particular, we use $\| \cdot\|$ instead of $\| \cdot\|_2$. 
\[f \in L^{1,\gamma}({\bf R}^{n}) \Leftrightarrow f \in L^{1}({\bf R}^{n}), \Vert f\Vert_{1,\gamma} := \int_{{\bf R}^{n}}(1+\vert x\vert)^{\gamma}\vert f(x)\vert dx < +\infty, \quad \gamma \geq 0.\]\\
Furthermore, we denote the Fourier transform $\hat{\phi}(\xi)$ of the function $\phi(x)$ by
\begin{equation}
{\cal F}(\phi)(\xi) := \hat{\phi}(\xi) := \frac{1}{(2\pi)^{n/2}}\int_{{\bf R}^{n}}e^{-ix\cdot\xi}\phi(x)dx,
\end{equation}
where $i := \sqrt{-1}$, and $x\cdot\xi = \displaystyle{\sum_{i=1}^{n}}x_{i}\xi_{i}$ for $x = (x_{1},\cdots,x_{n})$ and $\xi = (\xi_{1},\cdots,\xi_{n})$, and the inverse Fourier transform of ${\cal F}$ is denoted by ${\cal F}^{-1}$. When we estimate several functions by applying the Fourier transform sometimes we can also use the following definition in place of (1.5)
\[{\cal F}(\phi)(\xi) := \int_{{\bf R}^{n}}e^{-ix\cdot\xi}\phi(x)dx\]
without loss of generality. We also use the notation
\[v_{t}=\frac{\partial v}{\partial t}, \quad v_{tt}=\frac{\partial^{2} v}{\partial t^{2}}, \quad \Delta = \sum^n_{i=1}\frac{\partial^2}{\partial x_i^2},\ \ x=(x_1,\cdots,x_n),\]
and $\Vert D^{\ell}f\Vert := \left(\displaystyle{\int_{{\bf R}^{n}}}\vert\xi\vert^{2\ell}\vert\hat{f}(\xi)\vert^{2}\right)^{1/2}$ for $f \in H^{\ell}({\bf R}^{n})$.}


\section{Proofs of Theorems 1.1 and 1.2.}

To begin with, let us start with proving Theorem 1.1\\
In order to use the energy method in the Fourier space, we shall prepare the following notation.
\[{E_0(t,\xi):=\frac{1}{2}\vert \hat{u}_{t}\vert^2+\frac{1}{2}\vert\xi\vert^2\vert\hat{u}\vert^2},\]
\[E(t,\xi):=\frac{1}{2}\vert \hat{u}_{t}\vert^2+\frac{1}{2}\vert\xi\vert^2\vert\hat{u}\vert^2+\beta\rho(\xi)\mathfrak{R} (\hat{u}_t\overline{\hat{u}})+\frac{1}{2}\beta\rho(\xi)\vert\xi\vert^{2\theta}\vert\hat{u}\vert^2,\]
\[F(t,\xi):=\vert\xi\vert^{2\theta}\vert\hat{u}_t\vert^2+\beta\rho(\xi)\vert\xi\vert^2\vert\hat{u}\vert^2,\]
\[R(t,\xi):=\beta\rho(\xi)\vert\hat{u}_t\vert^2.\]
For $\theta > 1$, we define a key function $\rho: {\bf R}_{\xi}^{n}\rightarrow\bf R$ by
\[\rho(\xi)=\frac{\vert\xi\vert^{2\theta}}{1+\vert\xi\vert^{4\theta-2}}.\]
The discovery of this key function is a crucial point in this section. Note that most part of the following computations below can be also applied to the case of $\theta > 0$, for simplicity we restrict them only to the case of $\theta > 1$. 

Now, let us apply the Fourier transform to the both sides of (1.1) together with the initial data (1.2). Then in the Fourier space ${\bf R}_{\xi}^{n}$ one has the reduced problem
\begin{equation}
\hat{u}_{tt}(t,\xi)+\vert\xi\vert^2\hat{u}(t,\xi)+\vert \xi\vert^{2\theta}\hat{u}_t(t,\xi)=0,\quad(t,\xi)\in(0,\infty)\times{\bf R}_{\xi}^{n},
\end{equation}
\begin{equation}
\hat{u}(0,\xi)=\hat{u}_0(\xi),\quad  \hat{u}_t(0,\xi)=\hat{u}_1(\xi),\quad\xi\in{\bf R}_{\xi}^{n}.
\end{equation}
Multiply both sides of (2.1) by $\overline{\hat{u}_t}$, and further $\beta \rho(\xi)\overline{\hat{u}}$. Then, by taking the real part of the resulting identities one has
\begin{equation}
\frac{d}{dt}E_0(t,\xi)+\vert\xi\vert^{2\theta}\vert\hat{u}_t\vert^2=0,
\end{equation}
\begin{equation}
\frac{d}{dt}\{\beta\rho(\xi)\mathfrak{R}(\hat{u}_t\overline{\hat{u}})+\frac{1}{2}\beta\vert\xi\vert^{2\theta}\rho(\xi)\vert\hat{u}\vert^2\}+\beta\vert\xi\vert^2\rho(\xi)\vert\hat{u}\vert^2=\beta\rho(\xi)\vert\hat{u}_t\vert^2.
\end{equation}
By adding (2.3) and (2.4), one has
\begin{equation}
\frac{d}{dt}E(t,\xi)+F(t,\xi)=R(t,\xi).
\end{equation}
We prove
\begin{lem}  For $\beta>0$, it is true that
\[R(t,\xi)\leq\beta F(t,\xi),\quad\xi\in{\bf R}_{\xi}^{n}.\]
\end{lem}
{\it Proof.}\,For all $\theta>0$, it is true that
\begin{align*}
R(t,\xi)&=\beta\rho(\xi)\vert\hat{u}_t\vert^2
=\beta\frac{\vert\xi\vert^{2\theta}}{1+\vert\xi\vert^{4\theta-2}}\vert\hat{u}_t\vert^2
\leq\beta\vert\xi\vert^{2\theta}\vert\hat{u}_t\vert^2
\leq\beta F(t,\xi),
\end{align*}
which is valid for any $\xi\in{\bf R}_{\xi}^{n}$. 
\hfill
$\Box$
\par 
\vspace{0.2cm}
\par
It follows from (2.5) and Lemma 2.1 that
\begin{equation}
\frac{d}{dt}E(t,\xi)+(1-\beta)F(t,\xi)\leq 0,
\end{equation}
provided that the parameters $\beta>0$ are small enough.
\begin{lem}
There is a constant $M_1>0$ depending on $\beta>0$ such that for all $\xi\in{\bf R}_{\xi}^{n}$ with $\xi\neq0$ it follows that
\[{\rm (i)}\quad \frac{\rho(\xi)+\beta\rho(\xi)^2\vert\xi\vert^{-1}}{2\vert\xi\vert^{2\theta}}\leq M_1,\]
\[{\rm (ii)}\quad\frac{1}{2\beta}+\frac{\rho(\xi)}{2\vert\xi\vert}+\frac{\rho(\xi)\vert\xi\vert^{2\theta-2}}{2}\leq M_1,\]
\end{lem}
{\it Proof.}\\
\noindent
(i)\,Since
\[\frac{\rho(\xi)}{2\vert\xi\vert^{2\theta}}=\frac{1}{2\vert\xi\vert^{2\theta}}\cdot\frac{\vert\xi\vert^{2\theta}}{1+\vert\xi\vert^{4\theta-2}}\leq\frac{1}{2},\]
and
\begin{align*}
\frac{\beta\rho(\xi)^2\vert\xi\vert^{-1}}{2\vert\xi\vert^{2\theta}}&=\frac{\beta}{2\vert\xi\vert^{2\theta+1}}\cdot\frac{\vert\xi\vert^{4\theta}}{\left(1+\vert\xi\vert^{4\theta-2}\right)^2}=\frac{\beta}{2}\cdot\frac{\vert\xi\vert^{2\theta-1}}{\left(1+\vert\xi\vert^{4\theta-2}\right)}\cdot\frac{1}{1+\displaystyle\vert\xi\vert^{4\theta-2}}\\
&\leq\frac{\beta}{2}\cdot\frac{1}{\vert\xi\vert^{2\theta-1}+\frac{1}{\displaystyle\vert\xi\vert^{2\theta-1}}}\leq\frac{\beta}{4},
\end{align*}
the statement is true.\\
\noindent
(ii)We have
\[\frac{\rho(\xi)}{2\vert\xi\vert}=\frac{1}{2\vert\xi\vert}\frac{\vert\xi\vert^{2\theta}}{1+\vert\xi\vert^{4\theta-2}}=\frac{1}{2}\cdot\frac{1}{\vert\xi\vert^{2\theta-1}+\frac{1}{\displaystyle\vert\xi\vert^{2\theta-1}}}\leq\frac{1}{4},\]
and
\[\frac{\rho(\xi)\vert\xi\vert^{2\theta-2}}{2}=\frac{\vert\xi\vert^{2\theta-2}}{2}\frac{\vert\xi\vert^{2\theta}}{1+\vert\xi\vert^{4\theta-2}}=\frac{\vert\xi\vert^{4\theta-2}}{2(1+\vert\xi\vert^{4\theta-2})}\leq\frac{1}{2},\]
which implies the desired estimate.
\hfill
$\Box$
\par 
\vspace{0.2cm}
\par
\begin{lem}
There is a constant $M_2>0$ such that for all $\xi\in{\bf R}_{\xi}^{n}$,
it follows that 
\[\rho(\xi)E(t,\xi)\leq M_2 F(t,\xi).\]
\end{lem}
{\it Proof.}\,If $\xi\neq0$, since one has
\[\mathfrak{R}(\hat{u}_t\overline{\hat{u}})\leq\frac{1}{2}\left(\frac{\vert\hat{u}_{t}\vert^2}{\vert\xi\vert}+\vert\xi\vert\vert\hat{u}\vert^2\right),\]
it follows that
\begin{align*}
\rho(\xi)E(t,\xi)&\leq\frac{\rho(\xi)}{2}\vert\hat{u}_{t}\vert^2+\frac{\rho(\xi)}{2}\vert\xi\vert^2\vert\hat{u}\vert^2+\frac{\beta\rho(\xi)^2}{2}\vert\xi\vert^{2\theta}\vert\hat{u}\vert^2+\frac{\beta\rho(\xi)^2}{2}\frac{\vert\hat{u}_{t}\vert^2}{\vert\xi\vert}+\frac{\beta\rho(\xi)^2}{2}\vert\xi\vert\vert\hat{u}\vert^2\\
&\quad\leq(\frac{\rho(\xi)+\beta\rho(\xi)^2\vert\xi\vert^{-1}}{2\vert\xi\vert^{2\theta}})\vert\xi\vert^{2\theta}\vert\hat{u}_{t}\vert^2+(\frac{1}{2\beta}+\frac{\rho(\xi)}{2\vert\xi\vert}+\frac{\rho(\xi)\vert\xi\vert^{2\theta-2}}{2})\beta\rho(\xi)\vert\xi\vert^2\vert\hat{u}\vert^2,
\end{align*}
so that from Lemma 2.2 one has
\[\rho(\xi)E(t,\xi)\leq M_1 \vert\xi\vert^{2\theta}\vert\hat{u}_{t}\vert^2+M_1\beta\rho(\xi)\vert\xi\vert^2\vert\hat{u}\vert^2,\] 
which implies the desired estimate with $M_2:=M_1$ just defined in Lemma 2.2. 
According to the define of $E(t,\xi)$ and $F(t,\xi)$, above inequality also holds true with $\xi=0$.
\hfill
$\Box$
\par 
\vspace{0.2cm}
\par
Lemma 2.3 and (2.6) imply
\begin{equation}
\frac{d}{dt}E(t,\xi)+(1-\beta)\rho(\xi)M_2^{-1}E(t,\xi)\leq0
\end{equation}for any $\xi\in{\bf R}_{\xi}^{n}.$ From (2.7) we find 
\begin{equation}
E(t,\xi)\leq e^{-\alpha\rho(\xi)t}E(0,\xi),\quad\xi\in {\bf R}_{\xi}^{n},
\end{equation}
where $\alpha:=(1-\beta)M_2^{-1} > 0$ with small $\beta \in (0,1)$.

On the other hand, in the case of $\xi\neq0$, since we have 
\begin{equation}
\pm\beta\rho(\xi)\mathfrak{R}(\hat{u}_{t}\overline{\hat{u}})\leq\frac{\beta}{2}\vert\xi\vert^{2\theta}\rho(\xi)\vert\hat{u}\vert^2+\frac{\beta\rho(\xi)}{2}\frac{\vert\hat{u}_{t}\vert^2}{\vert\xi\vert^{2\theta}},
\end{equation}
it follows from the definition of $E(t,\xi)$ and (2.9) with minus sign that
\begin{equation}
E(t,\xi)\geq\frac{1}{2}\left(1-\displaystyle\frac{\beta\rho(\xi)}{\vert\xi\vert^{2\theta}}\right)\vert\hat{u}_{t}\vert^2+\frac{1}{2}\vert\xi\vert^2\vert\hat{u}\vert^2. 
\end{equation}
And also, we see that 
\begin{equation}
1-\frac{\beta\rho(\xi)}{\vert\xi\vert^{2\theta}}=\{1-\left(\frac{\beta}{\vert\xi\vert^{2\theta}}\cdot\frac{\vert\xi\vert^{2\theta}}{1+\vert\xi\vert^{4\theta-2}}\right)\}\geq(1-\beta)>0,
\end{equation}
for $\xi\neq0$ if we choose small $\beta\in (0,1)$. So, we obtain
\begin{equation}
E(t,\xi)\geq(1-\beta)E_0(t,\xi),\quad\xi\neq0.
\end{equation}
Since $E(t,0)=E_0(t,0),$ (2.12) holds true for all $\xi\in{\bf R}_{\xi}^{n}.$ Thus, from (2.8) one has
\begin{equation}
E_0(t,\xi)\leq(1-\beta)^{-1}e^{-\alpha\rho(\xi)t}E(0,\xi),\quad\xi\in{\bf R}_{\xi}^{n}.
\end{equation}
While, because of (2.9) with plus sign, for $\xi\neq0$ one has
\begin{equation}E(t,\xi)\leq\frac{1}{2}\vert\hat{u}_{t}\vert^2+\frac{1}{2}\vert\xi\vert^2\vert\hat{u}\vert^2+\beta\vert\xi\vert^{2\theta}\rho(\xi)\vert\hat{u}\vert^2+\frac{\beta}{2}\frac{\rho(\xi)}{\vert\xi\vert^{2\theta}}\vert\hat{u}_{t}\vert^2.
\end{equation}
Since $\displaystyle{\frac{\rho(\xi)}{\vert\xi\vert^{2\theta}}} \leq1, \xi\neq0,$ and 
\begin{align*}
\vert\xi\vert^{2\theta}\rho(\xi)\leq \vert\xi\vert^2,
\end{align*}
it follows that
\begin{equation}
E(t,\xi)\leq CE_0(t,\xi),\quad\xi\neq0,
\end{equation}
where $C>0$ is a constant depending on small $\beta$. The inequality (2.15) also holds true with $\xi=0$. By (2.13) and (2.15) with $t=0$ one has arrived at the significant estimate.
\begin{lem}
Let $\theta\in (1,\infty)$. Then, there is a constant $C=C(\beta)>$0 and $\alpha=\alpha(\beta)>0$ such that for all $\xi\in {\bf R}_{\xi}^{n}$ it is true that
\begin{align*}
E_0(t,\xi)\leq Ce^{-\alpha\rho(\xi)t}E_0(0,\xi).
\end{align*}
\end{lem}

{\it Proof of Theorem 1.1.}\,By lemma 2.4 and the Plancherel theorem one has
\begin{align}
\int_{{\bf R}_{x}^n}(\vert u_t\vert^2+\vert\nabla u\vert^2)dx \leq C\int_{{\bf R}_{\xi}^n}(\vert\hat{u}_{t}\vert^2+\vert\xi\vert^2\vert\hat{u}\vert^2)d\xi\leq C\int_{{\bf R}_{\xi}^n} E_0(t,\xi)d\xi\notag\\
\leq C(\int_{\vert\xi\vert\leq 1}+\int_{\vert\xi\vert\geq1})e^{-\alpha\rho(\xi)t}(\vert\hat{u}_{1}(\xi)\vert^2+\vert\xi\vert^2\vert\hat{u}_{0}(\xi)\vert^2)d\xi =: C(I_{low} + I_{high}).
\end{align}
We first prepare the following standard formula.
\begin{equation}
\int_{\vert\xi\vert\leq1}e^{-\alpha\vert\xi\vert^{2\theta}t}\vert\xi\vert^k d\xi\leq C(1+t)^{-\frac{k+n}{2\theta}}\quad(t \geq 0),
\end{equation}
for each $k \in {\bf N}\cup\{0\}$.

Now. let us start with estimating both $I_{low}$ and $I_{high}$ based on the shape of $\rho(\xi)$. In fact, since $1<\theta$, we see that
$$\vert\xi\vert\leq1\Rightarrow\frac{\vert\xi\vert^{2\theta}}{1+\vert\xi\vert^{4\theta-2}}\geq\frac{\vert\xi\vert^{2\theta}}{2}.$$
Therefore we see that
\[
I_{low}\leq\int_{\vert\xi\vert\leq1}e^{-\eta\vert\xi\vert^{2\theta}}\vert\hat{u}_{1}(\xi)\vert^2d\xi+\int_{\vert\xi\vert\leq1}e^{-\eta\vert\xi\vert^{2\theta}}\vert\xi\vert^2\vert\hat{u}_{0}\vert^2d\xi\]
\begin{equation}
\leq C\Vert u_1\Vert_{1}^2(1+t)^{-\frac{n}{2\theta}}+C\Vert u_0\Vert_1^2(1+t)^{-\frac{n+2}{2\theta}},
\end{equation}
where $\eta = \displaystyle{\frac{\alpha}{2}}$. Here, we have just used the formula (2.17).

Next, we shall estimate the high frequency part. This part is crucial in the case of $\theta > 1$. First, one has $$\vert\xi\vert\geq1\Rightarrow\frac{\vert\xi\vert^{2\theta}}{1+\vert\xi\vert^{4\theta-2}}\geq\frac{1}{2\vert\xi\vert^{2\theta-2}}.$$
Hence, it follows that
\begin{align*}
I_{high}&=\int_{\vert\xi\vert\geq1}e^{-\alpha\frac{\vert\xi\vert^{2\theta}}{1+\vert\xi\vert^{4\theta-2}}t}(\vert\hat{u}_{1}\vert^2+\vert\xi\vert^2\vert\hat{u}_{0}\vert^2)d\xi\leq\int_{\vert\xi\vert\geq1}e^{-\alpha\frac{t}{2\vert\xi\vert^{2\theta-2}}}(\vert\hat{u}_{1}\vert^2+\vert\xi\vert^2\vert\hat{u}_{0}\vert^2)d\xi\\
&\leq\left(\sup_{\vert\xi\vert\geq1}\frac{e^{-\alpha t/2\vert\xi\vert^{2\theta-2}}}{\vert\xi\vert^{2\ell}}\right)\int_{\vert\xi\vert\geq1}(\vert\xi\vert^{2\ell}\vert\hat{u}_{1}\vert^2+\vert\xi\vert^{2\ell+2}\vert\hat{u}_{0}\vert^2) d\xi\\
&\leq C(1+t)^{-\frac{\ell}{\theta-1}}(\Vert D^{\ell}u_{1}\Vert^2+\Vert D^{\ell+1}u_{0}\Vert^2),
\end{align*}
where we have just used the following fact to get the desired decay rate: if we set $x=\displaystyle{\frac{\sqrt{\alpha t}}{\vert\xi\vert^{\theta-1}}}$, then one can compute as follows:
\[\sup_{\vert\xi\vert\geq1}\left(\frac{e^{-\alpha t/2\vert\xi\vert^{2\theta-2}}}{\vert\xi\vert^{2\ell}}\right) = \sup_{\vert\xi\vert\geq1}\left(\frac{e^{-\alpha (1+t)/2\vert\xi\vert^{2\theta-2}}}{\vert\xi\vert^{2\ell}}e^{\frac{\alpha}{2\vert\xi\vert^{2\theta-2}}}\right)\]

\[\leq \sup_{\vert\xi\vert\geq1}\left(\frac{e^{-\alpha (1+t)/2\vert\xi\vert^{2\theta-2}}}{\vert\xi\vert^{2\ell}}\right)e^{\alpha/2} = e^{\alpha/2}\sup_{r \geq 1}\left(\frac{e^{-\frac{\alpha(1+t)}{2r^{2\theta-2}}}}{r^{2\ell}}\right)\]
\[= e^{\alpha/2}\alpha^{-\frac{\ell}{\theta-1}}(1+t)^{-\frac{\ell}{\theta-1}}\sup\left\{\frac{\sigma^{2\ell/(\theta-1)}}{e^{\sigma^{2}/2}}:\,\sqrt{\alpha(1+t)} \geq \sigma \geq 0\right\} \]
\[
\leq e^{\alpha/2}\alpha^{-\frac{\ell}{\theta-1}}(1+t)^{-\frac{\ell}{\theta-1}}\sup_{\sigma \geq 0}\left( \frac{\sigma^{2\ell/(\theta-1)}}{e^{\sigma^{2}/2}}\right)
\]
\begin{equation}
\leq C(1+t)^{-\frac{\ell}{\theta-1}},
\end{equation}
where $C = C(\ell,\theta,\alpha) > 0$ is a constant. This completes the proof of Theorem 1.1.
\hfill
$\Box$
\par 
\vspace{0.2cm}
\par

{\it Proof of Theorem 1.2.}\,To begin with, from Lemma 2.4, if $\xi \ne 0$, then it is true that
\begin{equation}
\vert\hat{u}(t,\xi)\vert^{2} \leq Ce^{-\alpha\rho(\xi)t}(\frac{\vert\hat{u}_{1}(\xi)\vert^{2}}{\vert\xi\vert^{2}} + \vert\hat{u}_{0}(\xi)\vert^{2}).
\end{equation}
Furthermore, if $1<\theta$ and $\vert\xi\vert\leq1$, then we see that
\[\frac{\vert\xi\vert^{2\theta}}{1+\vert\xi\vert^{4\theta-2}}\geq\frac{\vert\xi\vert^{2\theta}}{2},\]
so that by integrating (2.20) over $\{\delta\leq\vert\xi\vert\leq1\}$ with small $\delta>0$, one gets 
\begin{align}
\int_{\delta\leq\vert\xi\vert\leq1}\vert\hat{u}(t,\xi)\vert^2&\leq C\int_{\delta\leq\vert\xi\vert\leq1}e^{-\eta t\vert\xi\vert^{2\theta}}\frac{\vert\hat{u}_{1}\vert^2}{\vert\xi\vert^2}d\xi+C\int_{\delta\leq\vert\xi\vert\leq1}e^{-\eta t\vert\xi\vert^{2\theta}}\vert\hat{u}_{0}\vert^2d\xi\notag\\
&\leq\Vert u_1\Vert_1^2\int_{\delta\leq\vert\xi\vert\leq1}e^{-\eta t\vert\xi\vert^{2\theta}}\vert\xi\vert^{-2}d\xi+C\Vert u_0\Vert_1^2\int_{\delta\leq\vert\xi\vert\leq1}e^{-\eta t\vert\xi\vert^{2\theta}}d\xi \notag\\
&\leq C(1+t)^{-\frac{n-2}{2\theta}}\Vert u_1\Vert_1^2+C(1+t)^{-\frac{n}{2\theta}}\Vert u_0\Vert_1^2,
\end{align}
where $\eta=\displaystyle{\frac{\alpha}{2}}$ and $C>0$ is independent of any $\delta > 0$. Here, we just have used the following formula in the case when $n -1 \geq k$:
\[\int_{\delta\leq\vert\xi\vert\leq1}e^{-\eta t\vert\xi\vert^{2\theta}}\vert\xi\vert^{-k}d\xi = \int_{\delta\leq\vert\xi\vert\leq1}e^{-\eta(1+ t)\vert\xi\vert^{2\theta}}\vert\xi\vert^{-k}e^{\eta\vert\xi\vert^{2\theta}}d\xi \leq Ce^{\eta}(1+ t)^{-\frac{n-k}{2\theta}}\int_0^{+\infty}e^{-z^{2\theta}}z^{n-k-1}dz\]
\begin{equation}
\leq C(1 + t)^{-\frac{n-k}{2\theta}}\int_0^{+\infty}e^{-z^{2\theta}}z^{n-k-1}dz,
\end{equation}
for $t \geq 0$, where $C>0$ does not depend on any $\delta > 0$. By letting $\delta\downarrow0$ in (2.21) one has the low frequency estimate:\\
\begin{equation}
\int_{\vert\xi\vert\leq1}\vert\hat{u}(t,\xi)\vert^2 \leq C(1+t)^{-\frac{n-2}{2\theta}}\Vert u_1\Vert_1^2+C(1+t)^{-\frac{n}{2\theta}}\Vert u_0\Vert_1^2.
\end{equation}
On the other hand, if $\vert\xi\vert\geq1$, since we have once more
\[\frac{\vert\xi\vert^{2\theta}}{1+\vert\xi\vert^{4\theta-2}}\geq\frac{1}{2\vert\xi\vert^{2\theta-2}},\]
it follows from (2.20) that 
\[\int_{\vert\xi\vert\geq1}\vert\hat{u}(t,\xi)\vert^2d\xi\leq C\int_{\vert\xi\vert\geq1}e^{-\alpha t/2\vert\xi\vert^{2\theta-2}}\frac{\vert\hat{u}_{1}\vert^2}{\vert\xi\vert^2}d\xi+C\int_{\vert\xi\vert\geq1}e^{-\alpha t/2\vert\xi\vert^{2\theta-2}}\vert\hat{u}_{0}\vert^2d\xi\]
\[\leq C\left(\sup_{\vert\xi\vert\geq1}\frac{e^{-\alpha t/2\vert\xi\vert^{2\theta-2}}}{\vert\xi\vert^{2\ell}}\right)\int_{\vert\xi\vert\geq1}(\vert\xi\vert^{2\ell-2}\vert\hat{u}_{1}\vert^2+\vert\xi\vert^{2\ell}\vert\hat{u}_{0}\vert^2) d\xi\]
\begin{equation}
\leq C(1+t)^{-\frac{\ell}{\theta-1}}(\Vert D^{\ell-1}u_{1}\Vert^2+\Vert D^{\ell}u_{0}\Vert^2),
\end{equation}
where we have used (2.19) again. (2.23) and (2.24) imply the desired estimate.
\hfill
$\Box$
\par

\section{Proof of Theorems 1.3.}

In this section, let us prove Theorem 1.3 by employing the method due to \cite{Ik-3,Ik-4}. We first establish the asymptotic profile of the solutions in the low frequency region $\vert\xi\vert \ll 1$, which is essential ingredient. It should be emphasized that the lemma below holds true for all $n \geq 1$ (cf. Theorem 1.2).
\begin{lem}
Let $n\geq1$, and $\theta = 2$. Then, it is true that there exists a constant $C>0$ such that for $t \gg 1$
\begin{align*}
\int_{\vert\xi\vert\leq\delta_0}&\vert\mathcal{F}(u(t,\cdot))(\xi)-\{P_1e^{-t\vert\xi\vert^4/2}\frac{\sin{(t\vert\xi\vert)}}{\vert\xi\vert}+P_0e^{-t\vert\xi\vert^4/2}\cos{(t\vert\xi\vert)}\}\vert^2d\xi\\
&\leq C\Vert u_1\Vert_{1,1}^2 t^{-\frac{n}{4}}+C\Vert u_0\Vert_{1,1}^2t^{-\frac{n+2}{4}} 
\end{align*}
with small positive $\delta_0\ll1$.
\end{lem}

In order to prove Lemma 3.1 we apply the Fourier transform with respect to the space variable $x$ of the both sides of (1.1)-(1.2). Then in the Fourier space ${\bf R}_{\xi}^n$ one has the reduced problem:
\begin{align}
\hat{u}_{tt}(t,\xi)+\vert\xi\vert^2\hat{u}(t,\xi)+\vert \xi\vert^4\hat{u}_t(t,\xi)=0,\quad(t,\xi)\in(0,\infty)\times{\bf R}_{\xi}^n,\\ 
\hat{u}(0,\xi)=\hat{u}_0(\xi),\qquad\hat{u}_t(0,\xi)=\hat {u}_1(\xi),\quad \xi\in{\bf R}_{\xi}^n.
\end{align}
Let us solve (3.1)-(3.2) directly under the condition that $0<\vert\xi\vert\leq\delta_0\ll1$. In this case we get
\begin{align}
\hat{u}(t,\xi)&=\frac{\hat{u}_1(\xi)-\sigma_2\hat{u}_0(\xi)}{\sigma_1-\sigma_2}e^{\sigma_1t}+\frac{\hat{u}_0(\xi)\sigma_1-\hat{u}_1(\xi)}{\sigma_1-\sigma_2}e^{\sigma_2t}\\
&=\frac{e^{\sigma_1t}-e^{\sigma_2t}}{\sigma_1-\sigma_2}\hat{u}_1(\xi)+\frac{\sigma_1e^{\sigma_2t}-\sigma_2e^{\sigma_1t}}{\sigma_1-\sigma_2}\hat{u}_0(\xi),
\end{align}
where $\sigma_j\in\bf {C}$ $(j=1,2)$ have forms:
\begin{align*}
\sigma_1=\frac{-\vert\xi\vert^4+i\vert\xi\vert\sqrt{4-\vert\xi\vert^6}}{2}, \quad \sigma_2=\frac{-\vert\xi\vert^4-i\vert\xi\vert\sqrt{4-\vert\xi\vert^6}}{2}.
\end{align*}
In this connection, the smallness of $\vert\xi\vert$ is assumed to guarantee $4 - \vert\xi\vert^{6} > 0$. 

Now, we use the decomposition of the initial data based on the idea due to \cite{Ik-4}:
\begin{equation}
\hat{u}_j(\xi)=A_j(\xi)-iB_j(\xi)+P_j\quad(j=0,1),
\end{equation}
where
\[A_j(\xi):=\int_{{\bf R}^n}(\cos(x\cdot\xi)-1)u_j(x)dx,\qquad B_j(x):=\int_{{\bf R}^n}\sin(x\cdot\xi)u_j(x)dx,\quad(j=0,1).\]
From  (3.4) and (3.5) we see that 
\begin{align}
\hat{u}(t,\xi)=&P_1\left(\frac{e^{\sigma_1t}-e^{\sigma_2t}}{\sigma_1-\sigma_2}\right)+P_0\left(\frac{\sigma_1e^{\sigma_2t}-\sigma_2e^{\sigma_1t}}{\sigma_1-\sigma_2}\right)\notag\\
&+(A_1(\xi)-iB_1(\xi))\left(\frac{e^{\sigma_1t}-e^{\sigma_2t}}{\sigma_1-\sigma_2}\right)\notag\\
&+(A_0(\xi)-iB_0(\xi))\left(\frac{\sigma_1e^{\sigma_2t}-\sigma_2e^{\sigma_1t}}{\sigma_1-\sigma_2}\right),
\end{align}
for all $\xi$ satisfying $0<\vert\xi\vert\leq\delta_0.$
It is easy to check that
\begin{equation}
\frac{e^{\sigma_1t}-e^{\sigma_2t}}{\sigma_1-\sigma_2}=\frac{2e^{-t\vert\xi\vert^4/2}\sin(\frac{t\vert\xi\vert\sqrt{4-\vert\xi\vert^6}}{2})}{\vert\xi\vert\sqrt{4-\vert\xi\vert^6}}, 
\end{equation}
\begin{equation}
\frac{\sigma_1e^{\sigma_2t}-\sigma_2e^{\sigma_1t}}{\sigma_1-\sigma_2}=\frac{\vert\xi\vert^3e^{-t\vert\xi\vert^4/2}\sin(\frac{t\vert\xi\vert\sqrt{4-\vert\xi\vert^6}}{2})}{\sqrt{4-\vert\xi\vert^6}}+e^{-t\vert\xi\vert^4/2}\cos\left(\frac{t\vert\xi\vert\sqrt{4-\vert\xi\vert^6}}{2}\right).
\end{equation}
If we set
\[K_1(t,\xi):=P_0\frac{\vert\xi\vert^3e^{-t\vert\xi\vert^4/2}\sin(\frac{t\vert\xi\vert\sqrt{4-\vert\xi\vert^6}}{2})}{\sqrt{4-\vert\xi\vert^6}},\]
\[K_2(t,\xi):=(A_1(\xi)-iB_1(\xi))\left(\frac{e^{\sigma_1t}-e^{\sigma_2t}}{\sigma_1-\sigma_2}\right),\]
\[K_3(t,\xi):=(A_0(\xi)-iB_0(\xi))\left(\frac{\sigma_1e^{\sigma_2t}-\sigma_2e^{\sigma_1t}}{\sigma_1-\sigma_2}\right),\]
then it follows from (3.6), (3.7) and (3.8) that
\begin{align}
\hat{u}(t,\xi)=2&P_1\frac{e^{-t\vert\xi\vert^4/2}\sin(\frac{t\vert\xi\vert\sqrt{4-\vert\xi\vert^6}}{2})}{\vert\xi\vert\sqrt{4-\vert\xi\vert^6}}+P_0e^{-t\vert\xi\vert^4/2}\cos\left(\frac{t\vert\xi\vert\sqrt{4-\vert\xi\vert^6}}{2}\right)\notag\\
&+K_1(t,\xi)+K_2(t,\xi)+K_3(t,\xi),\quad0<\vert\xi\vert\leq\delta_0.
\end{align}
Let us shave the needless factors of the right hand side of (3.9) to get a precise shape of the asymptotic profile by using the mean value theorem. In fact, from the mean value theorem it follows that
\[2\frac{\sin(\frac{t\vert\xi\vert\sqrt{4-\vert\xi\vert^6}}{2})}{\vert\xi\vert\sqrt{4-\vert\xi\vert^6}}=\frac{2}{\sqrt{4-\vert\xi\vert^6}}\frac{\sin{(t\vert\xi\vert)}}{\vert\xi\vert}+t\left(\frac{\sqrt{4-\vert\xi\vert^6}-2}{\sqrt{4-\vert\xi\vert^6}}\right)\cos{(\epsilon(t,\xi)),}\]
\[\cos{\left(\frac{t\vert\xi\vert\sqrt{4-\vert\xi\vert^6}}{2}\right)}=\cos{(t\vert\xi\vert)}-t\vert\xi\vert\left(\frac{\sqrt{4-\vert\xi\vert^6}-2}{2}\right)\sin{(\eta(t,\xi))},\]
where
\[\epsilon(t,\xi):=\frac{t\vert\xi\vert\sqrt{4-\vert\xi\vert^6}}{2}\theta+t\vert\xi\vert(1-\theta),\quad\exists\theta\in(0,1),\]
\[\eta(t,\xi):=\frac{t\vert\xi\vert\sqrt{4-\vert\xi\vert^6}}{2}\theta'+t\vert\xi\vert(1-\theta'),\quad\exists\theta'\in(0,1),\] 
so that from (3.9) one has arrived at the identity
\begin{align}
\hat{u}(t,\xi)=2&P_1e^{-t\vert\xi\vert^4/2}\frac{1}{\sqrt{4-\vert\xi\vert^6}}\frac{\sin{(t\vert\xi\vert)}}{\vert\xi\vert}+P_0e^{-t\vert\xi\vert^4/2}\cos{(t\vert\xi\vert)}\notag\\
&+P_1e^{-t\vert\xi\vert^4/2}t\left(\frac{\sqrt{4-\vert\xi\vert^6}-2}{\sqrt{4-\vert\xi\vert^6}}\right)\cos{(\epsilon(t,\xi))}\notag\\
&-P_0e^{-t\vert\xi\vert^4/2}t\vert\xi\vert\left(\frac{\sqrt{4-\vert\xi\vert^6}-2}{2}\right)\sin{(\eta(t,\xi))}+\sum_{j=1}^3K_j(t,\xi).
\end{align}
Furthermore, we apply again the mean value theorem to get
\begin{equation}
\frac{2}{\sqrt{4-\vert\xi\vert^6}}=1+\frac{6(\theta'')^5\vert\xi\vert^6}{(4-(\theta'')^6\vert\xi\vert^6)^{\frac{3}{2}}}, \quad\exists\theta''\in(0,1), 
\end{equation}
so that from (3.10) and (3.11), in the case of $0<\vert\xi\vert\leq\delta_0$ we find that
\begin{align}
\hat{u}(t,\xi)=&P_1e^{-t\vert\xi\vert^4/2}\frac{\sin{(t\vert\xi\vert)}}{\vert\xi\vert}+P_0e^{-t\vert\xi\vert^4/2}\cos{(t\vert\xi\vert)}\notag\\
&+P_1K_4(t,\xi)-P_0K_5(t,\xi)+K_1(t,\xi)+K_2(t,\xi)\notag\\
&+K_3(t,\xi)+K_6(t,\xi), 
\end{align}
where
\[K_4(t,\xi):=e^{-t\vert\xi\vert^4/2}t\frac{\sqrt{4-\vert\xi\vert^6}-2}{\sqrt{4-\vert\xi\vert^6}}\cos{(\epsilon(t,\xi))},\]
\[K_5(t,\xi):=e^{-t\vert\xi\vert^4/2}t\vert\xi\vert\left(\frac{\sqrt{4-\vert\xi\vert^6}-2}{2}\right)\sin{(\eta(t,\xi))},\]
\[K_6(t,\xi)=P_1e^{-t\vert\xi\vert^4/2}\frac{\sin{(t\vert\xi\vert)}}{\vert\xi\vert}\frac{6(\theta'')^5\vert\xi\vert^6}{(4-(\theta'')^6\vert\xi\vert^6)^{\frac{3}{2}}}.\]\\

In the following part, let us check decay orders of the remainder terms $K_{j}(t,\xi)$ ($j = 1,2,3,4,5,6$) of (3.12) based on the formula (2.17).\\
\noindent
{\it \bf(I)\,Estimate for $ K_1(t,\xi).$}
\begin{align}
\int_{\vert\xi\vert\leq\delta_0}\vert K_1(t,\xi)\vert^2d\xi&\leq\vert P_0\vert^2\int_{\vert\xi\vert\leq\delta_0}\frac{\vert\xi\vert^6e^{-t\vert\xi\vert^4} \vert\sin{(\frac{t\vert\xi\vert\sqrt{4-\vert\xi\vert^6}}{2})}\vert^2}{4-\vert\xi\vert^6}d\xi\notag\\
&\leq\frac{\vert P_0\vert^2}{4-\delta_0^6}\int_{\vert\xi\vert\leq\delta_0}\vert\xi\vert^6e^{-t\vert\xi\vert^4}d\xi\leq\frac{\vert P_0\vert^2}{4-\delta_0^6}t^{-\frac{n+6}{4}}.
\end{align}

Next, we use the following property
\begin{equation}
\vert\sqrt{4-\vert\xi\vert^6}-2\vert=\frac{\vert\xi\vert^6}{2+\sqrt{4-\vert\xi\vert^6}}\leq\vert\xi\vert^6.
\end{equation}
{\it \bf(II)\,Estimate for $ K_4(t,\xi).$}\\
It follows from (3.14) that
\begin{align}
\int_{\vert\xi\vert\leq\delta_0}\vert K_4(t,\xi)\vert^2d\xi&\leq Ct^2\int_{\vert\xi\vert\leq\delta_0}e^{-t\vert\xi\vert^4}\frac{\vert\sqrt{4-\vert\xi\vert^6}-2\vert^2}{(4-\vert\xi\vert^6)}\vert\cos{(\epsilon(t,\xi))}\vert^2d\xi\notag\\
&\leq C\frac{t^2}{4-\delta_0^6}\int_{\vert\xi\vert\leq\delta_0}e^{-t\vert\xi\vert^4}\vert\xi\vert^{12}d\xi\leq\frac{C}{(4-\delta_0^6)}t^{-\frac{n}{4}-1}.
\end{align}
{\it \bf(III)\,Estimate for $ K_5(t,\xi).$}\\
Again it follows from (3.14) that
\begin{align}
\int_{\vert\xi\vert\leq\delta_0}\vert K_5(t,\xi)\vert^2d\xi&\leq Ct^2\int_{\vert\xi\vert\leq\delta_0}e^{-t\vert\xi\vert^4}\vert\xi\vert^2\left|\frac{\sqrt{4-\vert\xi\vert^6}-2}{2}\right|^2\vert\sin{(\eta(t,\xi))}\vert^2d\xi\notag\\
&\leq Ct^2\int_{\vert\xi\vert\leq\delta_0}e^{-t\vert\xi\vert^4}\vert\xi\vert^{14}d\xi\leq Ct^{-\frac{n+6}{4}}.
\end{align}
{\it \bf(IV)\,Estimate for $ K_6(t,\xi).$}\\
\begin{align}
\int_{\vert\xi\vert\leq\delta_0}\vert K_6(t,\xi)\vert^2d\xi&\leq\vert P_1\vert^2
\int_{\vert\xi\vert\leq\delta_0}e^{-t\vert\xi\vert^4}\frac{\sin^2{(t\vert\xi\vert)}}{\vert\xi\vert^2}\frac{36(\theta'')^{10}\vert\xi\vert^{12}}{(\vert4-(\theta'')^6\vert\xi\vert^6\vert^3)}d\xi\notag\\
&\leq C\frac{\vert P_1\vert^2}{\vert4-\delta_0^6\vert^3}\int_{\vert\xi\vert\leq\delta_0}e^{-t\vert\xi\vert^4}\vert\xi\vert^{10}d\xi\leq C\vert P_1\vert^2 t^{-\frac{n+10}{4}}.
\end{align}
In order to estimate further $K_j\hspace{1mm} (j=2,3)$, we prepare the following lemma, introduced in \cite{Ik-3}.
\begin{lem}
Let $n \geq1$. Then it holds that for all $\xi \in {\bf R}^{n}$
\[\vert A_j (\xi)\vert\leq L\vert\xi\vert\Vert u_j\Vert_{1,1}\quad(j=0,1),\]
\[\vert B_j (\xi)\vert\leq M\vert\xi\vert\Vert u_j\Vert_{1,1}\quad(j=0,1),\]
where
\[L:=\sup_{\theta\neq0}\frac{\vert1-\cos{\theta}\vert}{\vert\theta\vert}<+\infty,\quad M:=\sup_{\theta\neq0}\frac{\vert\sin{\theta}\vert}{\vert\theta\vert}<+\infty,\]
and both $A_j(\xi)$ and $B_j(\xi)$ are defined in {\rm (3.5)}.
\end{lem}
Basing on Lemma 3.2 we check decay rates of $K_{j}(t,\xi)$ with $j = 2,3$. This part is essential in this paper.\\
\noindent
{\it \bf(V)\,Estimate for $ K_2(t,\xi).$}\\
It follows from (3.7) and Lemma 3.2 that
\begin{align}
\int_{\vert\xi\vert\leq\delta_0}\vert K_2(t,\xi)\vert^2d\xi&\leq C(L^2+M^2)\Vert u_1\Vert_{1,1}^2\int_{\vert\xi\vert\leq\delta_0}\vert\xi\vert^2e^{-t\vert\xi\vert^4}\frac{\sin^2{(\frac{t\vert\xi\vert\sqrt{4-\vert\xi\vert^6}}{2})}}{\vert\xi\vert^2 (4-\vert\xi\vert^6)}d\xi,\notag\\
&\leq C\Vert u_1\Vert_{1,1}^2\int_{\vert\xi\vert\leq\delta_0}e^{-t\vert\xi\vert^4}d\xi\leq C\Vert u\Vert_{1,1}^2t^{-\frac{n}{4}}.
\end{align}
{\it \bf(VI)\,Estimate for $K_3(t,\xi)$.}\\
Similarly, 
\[\int_{\vert\xi\vert\leq\delta_0}\vert K_3(t,\xi)\vert^2d\xi \leq  C(L^2+M^2)\Vert u_0\Vert_{1,1}^2\int_{\vert\xi\vert\leq\delta_0}\vert\xi\vert^8e^{-t\vert\xi\vert^4}\frac{\sin^2(\frac{t\vert\xi\vert\sqrt{4-\vert\xi\vert^6}}{2})}{(4-\vert\xi\vert^6)}d\xi \]
\[+  C(L^2+M^2)\Vert u_0\Vert_{1,1}^2\int_{\vert\xi\vert\leq\delta_0}\vert\xi\vert^2e^{-t\vert\xi\vert^4}\cos^{2}\left(\frac{t\vert\xi\vert\sqrt{4-\vert\xi\vert^6}}{2}\right)d\xi\]
\[\leq C\Vert u_0\Vert_{1,1}^2\int_{\vert\xi\vert\leq\delta_0}\vert\xi\vert^8e^{-t\vert\xi\vert^4}d\xi+C\Vert u_0\Vert_{1,1}^2\int_{\vert\xi\vert\leq\delta_0}\vert\xi\vert^2e^{-t\vert\xi\vert^4}d\xi\]
\begin{equation}
\leq C\Vert u_0\Vert_{1,1}^2t^{-\frac{n}{4}-2}+C\Vert u_0\Vert_{1,1}^2 t^{-\frac{n+2}{4}}.
\end{equation}
Under these preparations let us prove Lemma 3.1. \\
\par
\vspace{0.2cm}
{\it Proof of Lemma 3.1.} It follows from (3.12), (3.13), (3.15), (3.16), (3.17), (3.18) and (3.19) that
\[\int_{\vert\xi\vert\leq\delta_0}\vert\mathcal{F}(u(t,\cdot))(\xi)-\{P_1e^{-t\vert\xi\vert^4/2}\frac{\sin{(t\vert\xi\vert)}}{\vert\xi\vert}+P_0e^{-t\vert\xi\vert^4/2}\cos{(t\vert\xi\vert)}\}\vert^2d\xi\]

\[\leq C\vert P_1\vert^2t^{-\frac{n}{4}-1}+C\vert P_1\vert^2t^{-\frac{n+10}{4}}+C\vert P_0\vert^2t^{-\frac{n+6}{4}}+\Vert u_1\Vert_{1,1}^2t^{-\frac{n}{4}}+C\Vert u_0\Vert_{1,1}^2t^{-\frac{n}{4}-2}+C\Vert u_0\Vert_{1,1}^2 t^{-\frac{n+2}{4}}\]
\[\leq C\Vert u_1\Vert_{1,1}^2 t^{-\frac{n}{4}}+C\Vert u_0\Vert_{1,1}^2t^{-\frac{n+2}{4}},\quad (t \gg 1)\]
which implies the desired estimate.
\hfill
$\Box$
\par 
\vspace{0.2cm}
\par

Finally, in this section let us prove our main Theorem 1.3 by relying on Lemma 2.4 and Lemma 3.1. The main part is the estimates in the high frequency region $\vert\xi\vert \gg 1$. In this part, we encounter the so called regularity loss type estimates, which are essentially a new point of view throughout this paper.\\
\noindent 
{\it Proof of Theorem 1.3.}
First we decompose the integrand to be estimated into two parts such that one is the low frequency part, and the other is high frequency one as follows.   
\begin{align*}
\int_{{\bf R}^n}&\vert\mathcal{F}(u(t,\cdot))(\xi)-\{P_{1}e^{-t\vert\xi\vert^4/2}\frac{\sin{(t\vert\xi\vert)}}{\vert\xi\vert}+P_0e^{-t\vert\xi\vert^4/2}\cos{(t\vert\xi\vert)}\}\vert^2d\xi\\
&=\left(\int_{\vert\xi\vert\leq\delta_0}+\int_{\vert\xi\vert\geq\delta_0}\right)\vert\mathcal{F}(u(t,\cdot))(\xi)-\{P_1e^{-t\vert\xi\vert^4/2}\frac{\sin{(t\vert\xi\vert)}}{\vert\xi\vert}+P_0e^{-t\vert\xi\vert^4/2}\cos{(t\vert\xi\vert)}\}\vert^2d\xi\\
&=:I_l(t)+I_h(t).
\end{align*}
As a direct consequence of Lemma 3.1 we can see soon that 
\begin{equation}
I_l(t)\leq C\Vert u_1\Vert_{1,1}^2 t^{-\frac{n}{4}}+C\Vert u_0\Vert_{1,1}^2t^{-\frac{n+2}{4}}.
\end{equation}

On the other hand, it follows from Lemma 2.4 with $\theta=2$ that
\[\vert\hat{u}_{t}(t,\xi)\vert^2+\vert\xi\vert^2\vert\hat{u}(t,\xi)\vert^2\leq Ce^{-\alpha\rho(\xi)t}(\vert\hat{u}_{1}(\xi)\vert^2+\vert\xi\vert^2\vert\hat{u}_{0}(\xi)\vert^2),\quad\xi\in{\bf R}^n.\]

In this case, we see that if $\delta_0\leq\vert\xi\vert\leq 1$, then 
$$\frac{\vert\xi\vert^{4}}{1+\vert\xi\vert^{6}}\geq\frac{\vert\xi\vert^{4}}{2}\geq\frac{\delta_0^4}{2},$$
and if $\vert\xi\vert\geq1$, then 
$$\frac{\vert\xi\vert^{4}}{1+\vert\xi\vert^{6}}\geq\frac{1}{2\vert\xi\vert^2}.$$
So, by similar argument to the proof of Theorem 1.2 (see (2.24)), which used (2.19) one can estimate as follows:
\[
\int_{\vert\xi\vert\geq\delta_0}\vert\hat{u}(t,\xi)\vert^2d\xi \leq C\int_{\vert\xi\vert\geq\delta_0}e^{-\alpha\rho(\xi)t}\left(\frac{\vert\hat{u}_{1}(\xi)\vert^2}{\vert\xi\vert^2}+\vert\hat{u}_{0}(\xi)\vert^2\right)d\xi\]
\[=\int_{1\geq\vert\xi\vert\geq\delta_0}e^{-\alpha\frac{\vert\xi\vert^{4}}{1+\vert\xi\vert^{6}}t}\left(\frac{\vert\hat{u}_{1}(\xi)\vert^2}{\vert\xi\vert^2}+\vert\hat{u}_{0}(\xi)\vert^2\right)d\xi +\int_{\vert\xi\vert\geq 1}e^{-\alpha\frac{\vert\xi\vert^{4}}{1+\vert\xi\vert^{6}}t}\left(\frac{\vert\hat{u}_{1}(\xi)\vert^2}{\vert\xi\vert^2}+\vert\hat{u}_{0}(\xi)\vert^2\right)d\xi\]
\begin{equation}
\leq Ce^{-c\alpha t}(\Vert u_1\Vert^2+\Vert u_0\Vert^2)
+Ct^{-\ell}(\Vert D^{\ell}u_{0}\Vert^2 + \Vert D^{\ell-1} u_{1}\Vert^2)\quad (t \gg 1),
\end{equation}
where the constants $C>0$ and $\alpha>0$ depend on $\delta_0.$ On the other hand, 
\[
P_1^2\int_{\vert\xi\vert\geq\delta_0}e^{-t\vert\xi\vert^4}\left|\frac{\sin{(t\vert\xi\vert)}}{\vert\xi\vert}\right|^2d\xi+P_0^2\int_{\vert\xi\vert\geq\delta_0}e^{-t\vert\xi\vert^4}\vert\cos{(t\vert\xi\vert)}\vert^2d\xi\]
\[\leq C\left(\frac{1}{\delta_0^2}+1\right)(\Vert u_1\Vert_1^2+\Vert u_0\Vert_1^2)\int_{\vert\xi\vert\geq\delta_0}e^{-t\vert\xi\vert^4}d\xi\]
\[=C\left(\frac{1}{\delta_0^2}+1\right)(\Vert u_1\Vert_1^2+\Vert u_0\Vert_1^2)\int_{\vert\xi\vert\geq\delta_0}e^{-\frac{t\vert\xi\vert^4}{2}}e^{-\frac{t\vert\xi\vert^4}{2}}d\xi\]
\[\leq C\left(\frac{1}{\delta_0^2}+1\right)(\Vert u_1\Vert_1^2+\Vert u_0\Vert_1^2)e^{-\frac{t\delta_0^4}{2}}\int_{\vert\xi\vert\geq\delta_0}e^{-\frac{t\vert\xi\vert^4}{2}}d\xi\]
\begin{equation}
\leq C(\Vert u_1\Vert_1^2+\Vert u_0\Vert_1^2)e^{-\frac{t\delta_0^4}{2}}t^{-\frac{n}{4}}\leq C(\Vert u_1\Vert_1^2+\Vert u_0\Vert_1^2)e^{-\alpha t}\quad(t \gg 1),
\end{equation}
with some $C>0$ and $\alpha>0$. Therefore, by evaluating $I_h(t)$ based on (3.21) and (3.22), and combining it with (3.20) one can arrive at
\begin{align*}
\int_{{\bf R}^n}&\vert\mathcal{F}(u(t,\cdot))(\xi)-\{P_1e^{-t\vert\xi\vert^4/2}\frac{\sin{(t\vert\xi\vert)}}{\vert\xi\vert}+P_0 e^{-t\vert\xi\vert^4/2}\cos{(t\vert\xi\vert)}\}\vert^2d\xi\\
&\leq C\Vert u_1\Vert_{1,1}^2 t^{-\frac{n}{4}}+C\Vert u_0\Vert_{1,1}^2 t^{-\frac{n+2}{4}} + Ce^{-\alpha t}(\Vert u_1\Vert_1^2+\Vert u_0\Vert_1^2+\Vert u_1\Vert^2+\Vert u_0\Vert^2)\\
&\qquad+Ct^{-\ell}(\Vert D^{\ell}u_{0}\Vert^2 + \Vert D^{\ell-1} u_{1}\Vert^2)\quad(t\gg 1),
\end{align*}
where $\ell \geq 0$. This implies the desired estimate. In this connection, the additional restriction on $\ell \geq 0$ such as in Theorem 1.3 is added further based on Lemmas 4.1 and 4.2 stated in section 4 in order to justify the statement of the main result.
\hfill
$\Box$
\par 
\vspace{0.2cm}
\par

\section{Optimality of the decay estimates.}

In this section, based on the result of Theorem 1.3, we shall study the optimality of the decay rates of the solution obtained in Theorem 1.2 with $\theta = 2$. The whole idea comes from \cite{IO}, which dealt with the case $\theta = 1$ (no regularity loss type). \\

First, in the case of $n\geq3$, one easily gets
\begin{align}
\int_{{\bf R}^n}e^{-t\vert\xi\vert^4}\frac{\sin^2{(t\vert\xi\vert)}}{\vert\xi\vert^2}d\xi\leq\int_{{\bf R}^n}e^{-t\vert\xi\vert^4}\vert\xi\vert^{-2}d\xi\leq Ct^{-\frac{n-2}{4}}\quad(t\gg 1).
\end{align}
While, we also see that
\begin{align}
\int_{{\bf R}^n}e^{-t\vert\xi\vert^4}\frac{\sin^2{(t\vert\xi\vert)}}{\vert\xi\vert^2}d\xi=\frac{1}{2}\left(\int_{\vert\omega\vert=1}d\omega\right)t^{-\frac{n-2}{4}}(A_0-F_n(t)),\notag
\end{align}
where
\[ A_0:=\int_0^{\infty}e^{-x^4}x^{n-3}dx,\quad  F_n(t):=\int_0^{\infty}e^{-x^4}x^{n-3}\cos(2t^{3/4}x)dx.\]
So, in the case of $n \geq 3$,  since $e^{-x^4}x^{n-3}\in L^1(0,\infty)$, it follows from the Riemann-Lebesgue theorem that
\[\lim_{t \to \infty}F_{n}(t) = 0,\]
so that for sufficiently large $t$ one can get
\[\vert F_{n}(t)\vert = \left|\int_0^{\infty}e^{-x^4}x^{n-3}\cos(2t^{3/4}x)dx\right| \leq\frac{A_0}{2}.\]
Hence, for $t \gg 1$ it is true that 
\begin{align}
\int_{{\bf R}^n}e^{-t\vert\xi\vert^4}\frac{\sin^2{(t\vert\xi\vert)}}{\vert\xi\vert^2}d\xi\geq Ct^{-\frac{n-2}{4}},\quad (n \geq 3)
\end{align}
with some constant $C > 0$. 

Next, in the case of $n\geq1$, we have
\begin{align}
\int_{{\bf R}^n}e^{-t\vert\xi\vert^4}\cos^2{(t\vert\xi\vert)}d\xi
\leq\int_{{\bf R}^n}e^{-t\vert\xi\vert^4}d\xi\leq Ct^{-\frac{n}{4}}\quad(t\gg 1).
\end{align}
Furthermore, if we set
\[ B_0:=\int_0^{\infty}e^{-x^4}x^{n-1}dx,\quad  G_n(t):=\int_0^{\infty}e^{-x^4}x^{n-1}\cos(2t^{3/4}x)dx,\]
then one has
\[\int_{{\bf R}^n}e^{-t\vert\xi\vert^4}\cos^2{(t\vert\xi\vert)}d\xi = \frac{1}{2}\left(\int_{\vert\omega\vert=1}d\omega\right)t^{-\frac{n}{4}}(B_0 + G_n(t)),\]
so that by similar argument to (4.2), and by summarizing the arguments above one can get
\begin{lem}\,There exist constants $0 < C_{1} \leq C_{2} < +\infty$ such that for $t \gg 1$
\[C_{1}t^{-\frac{n-2}{4}} \leq \int_{{\bf R}^n}e^{-t\vert\xi\vert^4}\frac{\sin^2{(t\vert\xi\vert)}}{\vert\xi\vert^2}d\xi \leq C_{2}t^{-\frac{n-2}{4}}, \quad (n \geq 3),\]
\[C_{1}t^{-\frac{n}{4}} \leq \int_{{\bf R}^n}e^{-t\vert\xi\vert^4}\cos^2{(t\vert\xi\vert)}d\xi \leq C_{2}t^{-\frac{n}{4}},\quad (n \geq 1).\]
\end{lem}

In order to obtain the optimal estimate of the solution in the low dimensional case, i.e., $n = 1,2$, we prove the following lemma.
\begin{lem}
There exist constants $0<C_1\leq C_2\leq\infty$ such that, for sufficiently large $t>0$, it is true that
\begin{equation}
C_1t\leq\int_{-\infty}^{\infty}e^{-t\vert\xi\vert^4}\frac{\sin^2{(t\vert\xi\vert)}}{\vert\xi\vert^2}d\xi\leq C_2t,\quad (n = 1),
\end{equation}
\begin{equation}
C_1\log{t}\leq\int_{{\bf R}^2}e^{-t\vert\xi\vert^4}\frac{\sin^2{(t\vert\xi\vert)}}{\vert\xi\vert^2}d\xi\leq C_2\log{t},\quad (n = 2).
\end{equation}
\end{lem}

Before going to the proof, let us point out that the change of variables $\eta=t^{\frac{1}{4}}\xi$ yields that
\begin{equation}
\int_{{\bf R}^n}e^{-t\vert\xi\vert^4}\frac{\sin^2{(t\vert\xi\vert)}}{\vert\xi\vert^2}d\xi=t^{-\frac{n-2}{4}}\int_{{\bf R}^n}e^{-\vert\eta\vert^4}\frac{\sin^2{(t^{{3/4}}\vert\eta\vert)}}{\vert\eta\vert^2}d\eta=:t^{-\frac{n-2}{4}}I^{(n)}(t),
\end{equation}
for $n\geq1.$\\

\noindent
{\it Proof of lemma 4.2.}\,Let us divide the integral $I^{(n)}$ as
\[I^{(n)}(t) =\int_{\{\vert\eta\vert\leq t^{-\frac{3}{4}}\}}e^{-\vert\eta\vert^4}\frac{\sin^2{(t^{{3/4}}\vert\eta\vert)}}{\vert\eta\vert^2}d\eta + \int_{\{\vert\eta\vert> t^{-\frac{3}{4}}\}}e^{-\vert\eta\vert^4}\frac{\sin^2{(t^{{3/4}}\vert\eta\vert)}}{\vert\eta\vert^2}d\eta\]
\begin{equation}
:=I_1^{(n)}(t)+I_2^{(n)}(t).
\end{equation}
We first treat the case $n = 1$.\\

For the estimate of $I_1^{(1)}(t)$ from above, by using $\sin{(t^{\frac{3}{4}}\eta)}\leq t^{\frac{3}{4}}\eta$ for $0\leq\eta\leq t^{-\frac{3}{4}}$ we obtain
\begin{align}
I_1^{(1)}(t)&\leq 2t^{\frac{3}{2}}\int_0^{t^{-\frac{3}{4}}}e^{-\eta^4}d\eta=2t^{\frac{3}{2}}\biggl\{\bigl[\eta e^{-\eta^4}\bigr]_0^{t^{-\frac{3}{4}}}
+\int_0^{t^{-\frac{3}{4}}}4\eta^4e^{-\eta^4}d\eta\biggr\}\notag\\
&\leq2t^{\frac{3}{2}}\left( t^{-\frac{3}{4}}e^{-t^{-3}}+\int_0^{t^{-\frac{3}{4}}}4 t^{-3}e^{-t^{-3}}d\eta\right)\leq t^{\frac{3}{4}}+O(t^{-\frac{9}{4}}),\,(t \to \infty),
\end{align}
where we just have used the monotone increasing property of the function $\eta^4e^{-\eta^4}$ on $0\leq\eta\leq t^{-\frac{3}{4}}$ for $t\geq1$.

Similarly, since $2\sin{(t^{\frac{3}{4}}\eta)}\geq t^{\frac{3}{4}}\eta$ for $0\leq\eta\leq t^{-\frac{3}{4}},$ we can estimate $I_1^{(1)}(t)$ from below as follows:
\begin{align}
I_1^{(1)}(t)&\geq\frac{t^{\frac{3}{2}}}{4}\int_0^{t^{-\frac{3}{4}}}e^{-\eta^4}d\eta=\frac{t^{\frac{3}{2}}}{4}\biggl\{\biggl[\eta e^{-\eta^4}\biggr]_0^{t^{-\frac{3}{4}}}
+\int_0^{t^{-\frac{3}{4}}}4\eta^4e^{-\eta^4}d\eta\biggr\}\notag\\
&\geq\frac{t^{\frac{3}{2}}}{4} t^{-\frac{3}{4}}e^{-t^{-3}} = \frac{t^{\frac{3}{4}}}{4} + O(t^{-\frac{9}{4}}),
\end{align}as $t\rightarrow\infty$.

On the other hand, the estimate of $I_2^{(1)}(t)$ from above is obtained as follows:
\begin{align}
I_2^{(1)}(t)&\leq2\int_{t^{-\frac{3}{4}}}^{\infty}\frac{e^{-\eta^4}}{\eta^2}d\eta=\biggl[-\frac{2e^{-\eta^4}}{\eta}\biggr]_{t^{-\frac{3}{4}}}^{\infty}-\int_{t^{-\frac{3}{4}}}^{\infty}8\eta^2e^{-\eta^4}d\eta\notag\\
&\leq2 t^{\frac{3}{4}}e^{-t^{-3}}=2t^{\frac{3}{4}}+O(t^{-\frac{9}{4}}),\quad (t \to \infty).
\end{align}
Furthermore, in order to estimate $I_2^{(1)}(t)$ from below, we observe that
\[\left|\sin{(t^{\frac{3}{4}}\eta)}\right|\geq\frac{1}{\sqrt{2}}\qquad (\frac{1}{4}+j)\frac{\pi}{t^{\frac{3}{4}}}\leq\eta\leq(\frac{3}{4}+j)\frac{\pi}{t^{\frac{3}{4}}},\, (j=1,2,\cdots). \]
Now, since $\displaystyle{\frac{e^{-\eta^4}}{\eta^2}}$ is monotone decreasing, we obtain
\[I_2^{(1)}(t)\geq\sum_{j=1}^{\infty}\int_{(\frac{1}{4}+j)\frac{\pi}{t^{3/4}}}^{(\frac{3}{4}+j)\frac{\pi}{t^{3/4}}}\frac{e^{-\eta^4}}{\eta^2}d\eta\geq\frac{1}{2}\int_{\frac{5\pi}{4t^{3/4}}}^{\infty}\frac{e^{-\eta^4}}{\eta^2}d\eta = \frac{1}{2}\biggl[-\frac{e^{-\eta^4}}{\eta}\biggr]_{{\frac{5\pi}{4t^{3/4}}}}^{\infty}-2\int_{{\frac{5\pi}{4t^{3/4}}}}^{\infty}\eta^2e^{-\eta^4}d\eta\]
\begin{equation}
\geq\frac{2t^{3/4}}{5\pi}e^{-625\pi^4/256t^3}-2\int_0^{\infty}\eta^2e^{-\eta^4}d\eta=\frac{2t^{3/4}}{5\pi}+O(1),\quad (t \to \infty).
\end{equation}
Combining (4.6)- (4.11) yields (4.4).

Next, in the following paragraph we shall prove (4.5).

In fact, similarly to the one dimensional case, we obtain the estimate for $I_1^{(2)}(t)$ from above and below by using the relation
$$\displaystyle{\frac{1}{2}}t^{\frac{3}{4}}\vert\eta\vert\leq\sin{(t^{\frac{3}{4}}\vert\eta\vert)}\leq t^{\frac{3}{4}}\vert\eta\vert$$
for $0\leq\vert\eta\vert\leq t^{-\frac{3}{4}}$. The polar co-ordinate transform implies
\[I_1^{(2)}(t)\leq t^{\frac{3}{2}}\int_{\vert\eta\vert\leq t^{-\frac{3}{4}}}e^{-\vert\eta\vert^4}d\eta= 2\pi t^{\frac{3}{2}}\int_0^{t^{-\frac{3}{4}}}e^{-r^4}r dr = 2\pi t^{\frac{3}{2}}\biggl\{\biggl[\frac{1}{2}r^2e^{-r^4}\biggr]_0^{t^{-\frac{3}{4}}}+\int_0^{t^{-\frac{3}{4}}}2r^5e^{-r^4}dr\biggr\}\]
\begin{equation}
\leq \pi t^{\frac{3}{2}}t^{-\frac{3}{2}}e^{-t^{-3}}+4\pi t^{\frac{3}{2}}\int_0^{t^{-\frac{3}{4}}}t^{-15/4}e^{-t^{-3}}dr = \pi+O(t^{-3}),
\end{equation}
where we have used again the monotone increasing property of the function $r^5e^{-r^4}$ on $0\leq r\leq t^{-\frac{3}{4}}$ with large $t$.
Furthermore, we see that
\[I_1^{(2)}(t)\geq \frac{t^{\frac{3}{2}}}{4}\int_{\vert\eta\vert\leq t^{-\frac{3}{4}}}e^{-\vert\eta\vert^4}d\eta=\frac{\pi}{2}t^{\frac{3}{2}}\int_0^{t^{-\frac{3}{4}}}e^{-r^4}r dr = \frac{\pi}{2}t^{\frac{3}{2}}\biggl\{\biggl[\frac{1}{2}r^2e^{-r^4}\biggr]_0^{t^{-\frac{3}{4}}}+\int_0^{t^{-\frac{3}{4}}}2r^5e^{-r^4}dr\biggr\}\]
\begin{equation}
\geq\frac{\pi}{4}t^{\frac{3}{2}}t^{-\frac{3}{2}}e^{-t^{-3}} = \frac{\pi}{4}e^{-t^{-3}} \geq \frac{\pi}{8},\quad (t \gg 1).
\end{equation}

Finally, we can obtain the estimate for $I_2^{(2)}(t)$ from above and below by using the polar co-ordinate transform again:
\[I_2^{(2)}(t) \leq\int_{\vert\eta\vert\geq t^{-\frac{3}{4}}}\frac{e^{-\vert\eta\vert^4}}{\vert\eta\vert^2}d\eta=2\pi\int_{t^{-\frac{3}{4}}}^{\infty}\frac{e^{-r^4}}{r}dr = 2\pi\left(\biggl[e^{-r^4}\log{r}\biggr]_{t^{-\frac{3}{4}}}^{\infty}+\int_{t^{-\frac{3}{4}}}^{\infty}4r^3e^{-r^4}\log{r}dr\right)\]
\begin{equation}
\leq2\pi\log{t^{\frac{3}{4}}}e^{-t^{-3}}+8\pi\int_{0}^{\infty}r^3e^{-r^4}\vert\log{r}\vert dr \leq 3\pi\log{t}+O(1),\quad (t \gg 1).
\end{equation}
The same argument can be applied to the estimate of $I_2^{(2)}(t)$ from below as follows:
\[I_2^{(2)}(t) \geq \sum_{j=1}^{\infty}\frac{1}{2}\int_{(\frac{1}{4}+j)\frac{\pi}{t^{\frac{3}{4}}} \leq\vert\eta\vert\leq(\frac{3}{4}+j)\frac{\pi}{t^{\frac{3}{4}}}}\frac{e^{-\vert\eta\vert^4}}{\vert\eta\vert^2}d\eta
=\pi\sum_{j=1}^{\infty}\int_{(\frac{1}{4}+j)\frac{\pi}{t^{\frac{3}{4}}}}^{(\frac{3}{4}+j)\frac{\pi}{t^{\frac{3}{4}}}}\frac{e^{-r^4}}{r}dr \geq \frac{\pi}{2}\int_{\frac{5\pi}{4}t^{-\frac{3}{4}}}^{\infty}\frac{e^{-r^4}}{r}dr\]
\[=\frac{\pi}{2}\left(\biggl[e^{-r^4}\log{r}\biggr]_{\frac{5\pi}{4}t^{-\frac{3}{4}}}^{\infty}+\int_{\frac{5\pi}{4}t^{-\frac{3}{4}}}^{\infty}4r^3e^{-r^4}\log{r}dr\right)\]
\begin{equation}
\geq\frac{\pi}{2}\left(\log{t^{\frac{3}{4}}}-\log{\frac{5\pi}{4}}\right)-2\pi\int_0^{\infty}r^3e^{-r^4}\vert\log{r}\vert dr = \frac{3\pi}{8}\log{t}+O(1),
\end{equation}
as $t\rightarrow\infty$. (4.12)-(4.15) imply the validity of (4.5).
\hfill
$\Box$
\par 
\vspace{0.2cm}
\par
Now, recall a series of inequalities below, which is based on the Plancherel theorem:
\begin{align}
C\Vert u(t,\cdot)\Vert&\geq\Vert\mathcal{F}(u(t,\cdot))(\xi)\Vert\geq\Vert P_1e^{-t\vert\xi\vert^4/2}\frac{\sin{(t\vert\xi\vert)}}{\vert\xi\vert}+P_0e^{-t\vert\xi\vert^4/2}\cos{(t\vert\xi\vert)}\Vert\notag\\
&\quad -\Vert\mathcal{F}(u(t,\cdot))(\xi)-\{P_1e^{-t\vert\xi\vert^4/2}\frac{\sin{(t\vert\xi\vert)}}{\vert\xi\vert}+P_0e^{-t\vert\xi\vert^4/2}\cos{(t\vert\xi\vert)\}}\Vert\notag\\
&\geq\vert P_1\vert\Vert e^{-t\vert\xi\vert^4/2}\frac{\sin{(t\vert\xi\vert)}}{\vert\xi\vert}\Vert-\vert P_0\vert\Vert e^{-t\vert\xi\vert^4/2}\cos{(t\vert\xi\vert)}\Vert\notag\\
&\quad-\Vert\mathcal{F}(u(t,\cdot))(\xi)-\{P_1e^{-t\vert\xi\vert^4/2}\frac{\sin{(t\vert\xi\vert)}}{\vert\xi\vert}+P_0e^{-t\vert\xi\vert^4/2}\cos{(t\vert\xi\vert)\}}\Vert,
\end{align}
and 
\[C\Vert u(t,\cdot)\Vert \leq \Vert\hat{u}(t,\cdot)\Vert \leq \Vert\mathcal{F}(u(t,\cdot))(\xi)-\{P_1e^{-t\vert\xi\vert^4/2}\frac{\sin{(t\vert\xi\vert)}}{\vert\xi\vert}+P_0e^{-t\vert\xi\vert^4/2}\cos{(t\vert\xi\vert)\}}\Vert\]
\begin{equation}
+ \vert P_1\vert\Vert e^{-t\vert\xi\vert^4/2}\frac{\sin{(t\vert\xi\vert)}}{\vert\xi\vert}\Vert + \vert P_0\vert\Vert e^{-t\vert\xi\vert^4/2}\cos{(t\vert\xi\vert)}\Vert.
\end{equation}

Then, basing on (4.16) and (4.17), by the use of Theorem 1.3 and Lemmas 4.1 and 4.2, we can get the following results, which show the optimality of the $L^{2}$-decay rates of the solution obtained in Theorem 1.2 for $n \geq 3$. In particular, the case of $n = 1,2$ implies the blow-up at infinite time. The argument is almost similar to that of \cite{Ik-4,IO} except for the regularity assumed on the initial data. To state our final results we set
\[I_{0} := (\Vert u_{0}\Vert + \Vert u_{0}\Vert_{1} + \Vert u_{0}\Vert_{1,1} + \Vert D^{\ell}u_{0}\Vert + \Vert u_{1}\Vert + \Vert u_{1}\Vert_{1} + \Vert u_{1}\Vert_{1,1} + \Vert D^{\ell-1}u_{1}\Vert).\]
\begin{theo}
Let $n\geq3$. Under the same assumptions as in Theorem {\rm 1.3}, if $[u_0,u_1]\in(H^{\ell}({\bf R}^n) \cap L^{1,1}({\bf R}^n))\times(H^{\ell-1}({\bf R}^n)\cap L^{1,1}({\bf R}^n))$, then it is true that
\[C_1\vert P_1\vert t^{-\frac{n-2}{8}}\leq \Vert u(t,\cdot)\Vert\leq C_2 I_{0}t^{-\frac{n-2}{8}},\]
for $t\gg1$, where $C_{j} > 0$ {\rm (}$j = 1,2${\rm )} are some constants.
\end{theo}
\begin{theo}
Let $n=2$. Under the same assumptions as in Theorem {\rm 1.3}, if $[u_0,u_1]\in(H^{\ell}({\bf R}^2) \cap L^{1,1}({\bf R}^2))\times(H^{\ell-1}({\bf R}^2)\cap L^{1,1}({\bf R}^2))$, then it is true that
\[C_1\vert P_1\vert \sqrt{\log{t}}\leq \Vert u(t,\cdot)\Vert\leq C_2 I_{0}\sqrt{\log{t}},\]
for $t\gg1$, where $C_{j} > 0$ {\rm (}$j = 1,2${\rm )} are some constants.
\end{theo}
\begin{theo}
Let $n=1$. Under the same assumptions as in Theorem {\rm 1.3}, if $[u_0,u_1]\in(H^{\ell}({\bf R}) \cap L^{1,1}({\bf R}))\times(H^{\ell-1}({\bf R})\cap L^{1,1}({\bf R}))$, then it is true that
\[C_1\vert P_1\vert \sqrt{t}\leq \Vert u(t,\cdot)\Vert\leq C_2 I_{0}\sqrt{t},\]
for $t\gg1$, where $C_{j} > 0$ {\rm (}$j = 1,2${\rm )} are some constants.
\end{theo}
\begin{rem}{\rm The big difference between the known results for $\theta = 1$ and this case of $\theta = 2$ is in the regularity assumed on the initial data to get the similar optimality. This discovery is the core of this paper.}
\end{rem}
\par
\vspace{0.2cm}

\noindent{\em Acknowledgement.}
\smallskip
The work of the first author (R. IKEHATA) was supported in part by Grant-in-Aid for Scientific Research (C) 15K04958 of JSPS.


\end{document}